\documentclass[11pt]{article}
\usepackage{amsmath, amssymb, amsthm, mathrsfs, graphicx, tikz, enumitem, hyperref, geometry}
\newtheorem{theorem}{Theorem}

\title{Generalizations of Dini’s Theorem under Weakened Monotonicity Conditions}
\author{Riwaj Khatiwada}
\date{}

\begin{document}

\maketitle

\begin{abstract}
Dini’s Theorem classically guarantees that a monotone sequence of continuous functions
converging pointwise on a compact interval to a continuous limit in fact converges uniformly. In this paper,
we establish new theorems generalizing Dini’s result by replacing the restrictive monotonicity assumption
with more flexible conditions – notably equicontinuity, convexity, or controlled variation hypotheses. We
prove that under these alternative sufficient conditions, pointwise convergence of continuous functions on
a compact interval still implies uniform convergence. Rigorous definitions and statements are provided,
followed by formal proofs. We construct illustrative examples and counterexamples to demonstrate when
pointwise limits do or do not become uniform limits under various assumptions. Graphical visualizations
accompany these examples, highlighting convergence behaviors. Our results situate within classical real
analysis – for context we review Dini’s Theorem as found in standard references (e.g. Rudin’s Principles of
Mathematical Analysis and Apostol’s Mathematical Analysis) – and extend the theory of uniform convergence
beyond monotone sequences.
\end{abstract}

\section{Introduction}

One of the subtle distinctions in real analysis is between \textbf{pointwise} and \textbf{uniform convergence} of function sequences. Recall that a sequence of functions $(f_n)$ on a domain $D$ converges \emph{pointwise} to a function $f$ if for every $x \in D$, we have
\[
f_n(x) \to f(x) \quad \text{as } n \to \infty.
\]
By contrast, \emph{uniform convergence} requires convergence to $f$ at a uniform rate over the whole domain: $(f_n)$ converges uniformly to $f$ on $D$ if for every $\varepsilon > 0$, there exists $N$ such that
\[
|f_n(x) - f(x)| < \varepsilon \quad \text{for all } x \in D \text{ whenever } n \ge N.
\]

Uniform convergence is a strong condition that, among other things, preserves continuity: if each $f_n$ is continuous and $f_n \to f$ uniformly, then the limit $f$ is also continuous. Pointwise convergence alone is generally too weak to ensure such preservation of continuity or interchange of limiting processes (integration, differentiation, etc.).

A classical example is the sequence $f_n(x) = x^n$ on $[0,1]$. For each fixed $x \in [0,1)$, we have $x^n \to 0$ as $n \to \infty$, and at $x = 1$ we have $1^n = 1$ for all $n$. So $f_n(x)$ converges pointwise to the limit function
\[
f(x) = 
\begin{cases}
0, & 0 \le x < 1, \\
1, & x = 1,
\end{cases}
\]
which is discontinuous on $[0,1]$. In fact, the convergence is not uniform in this case — if it were uniform, the limit would necessarily be continuous. This example shows that without additional conditions, pointwise convergence need not be uniform.

\medskip

\textbf{Dini’s Theorem} provides one powerful additional condition: \emph{monotonicity in $n$}. In its standard form, Dini’s Theorem states that if $K$ is a compact space (such as a closed interval) and $(f_n)$ is a monotone sequence of continuous real-valued functions on $K$ that converges pointwise to a continuous function $f$, then $f_n$ actually converges to $f$ uniformly on $K$. 

(The monotonicity can be either non-decreasing or non-increasing; in each case, one assumes $f$ is continuous on $K$, which, as noted, is necessary for uniform convergence.)

Informally, Dini’s Theorem highlights that greater control on how $f_n$ approaches $f$ — in this case, no oscillations since each $f_n(x)$ moves monotonically toward $f(x)$ for each $x$ — forces the convergence to be uniform. Indeed, as several textbooks remark, Dini’s scenario is “one of the very few situations where pointwise convergence (of continuous functions) implies uniform convergence”.

\subsection*{Motivation and Aim}
The monotonicity hypothesis in Dini’s Theorem, while useful, is quite restrictive. Many naturally arising function sequences are not monotonic in $n$ but still exhibit convergence with properties strong enough to enforce uniformity. The aim of this paper is to identify and prove alternative sufficient conditions—\emph{weaker or different than monotonicity}—that guarantee pointwise convergence implies uniform convergence. In particular, we explore three types of conditions:

\medskip

\textbf{1. Equicontinuity of the sequence:} Intuitively, equicontinuity means that all functions $f_n$ oscillate in $x$ at the same controlled rate. We will show that if $(f_n)$ is an equicontinuous sequence of continuous functions on a compact interval and $f_n$ converges pointwise to a continuous $f$, then $f_n$ converges uniformly to $f$. This result can be viewed as a direct consequence of the Arzelà–Ascoli Theorem and its corollaries, but we will also give an $\varepsilon$–$\delta$ style proof for completeness.

\medskip

\textbf{2. Convexity of the functions:} Convexity is a geometric condition on each $f_n$ (not on the sequence as $n$ changes). We will prove that if each $f_n$ is convex (and continuous) on a closed interval $[a,b]$, and $f_n \to f$ pointwise with $f$ continuous on $[a,b]$, then in fact $f_n \to f$ uniformly. Although the $f_n$ need not be monotone in $n$, convexity will be seen to provide a form of ``one-sided control'' on oscillations in $x$ that is sufficient to apply a Dini-like argument on small subintervals.

\medskip

\textbf{3. Structured bounded variation:} Functions of bounded variation can oscillate, but only in a limited total amount. We consider a uniformly bounded variation condition—roughly, assume there is a fixed $M$ such that each $f_n$ has total variation $\le M$ on $[a,b]$. This by itself does not guarantee uniform convergence (we will present a counterexample), so we refine the condition to require that the variation of $f_n$ cannot concentrate on arbitrarily small subintervals as $n$ grows. Under this uniformly distributed variation condition (made precise later), we will show pointwise convergence to a continuous limit again implies uniform convergence. This condition can be regarded as a quantitative relaxation of monotonicity: monotone functions have bounded variation (in fact equal to their range), and our condition permits a controlled number of oscillations of limited amplitude.

\medskip

By establishing the above results, we generalize Dini’s Theorem to broader contexts. These generalizations are not only of theoretical interest but also useful in analysis and applied settings (e.g., approximation of functions, convergence of algorithms) where one often wants to ensure uniform convergence without strict monotonicity. We situate our work in the context of classical real analysis—relying on standard theorems like Arzelà–Ascoli, properties of convex functions, and basic measure of variation (Jordan’s theorem)—and we cite references such as \textit{Rudin’s Principles of Mathematical Analysis} and \textit{Apostol’s Mathematical Analysis} for foundational definitions and results.

\medskip

The remainder of the paper is organized as follows. In Section 2, we review definitions and preliminary facts (uniform convergence, equicontinuity, convexity, bounded variation) for self-containedness. In Section 3, we formally state and prove the new theorems: \textbf{Theorem 1} (Equicontinuous Convergence Criterion), \textbf{Theorem 2} (Convex Sequence Criterion), and \textbf{Theorem 3} (Uniformly Distributed Bounded Variation Criterion). Section 4 provides examples illustrating each theorem and counterexamples showing necessity of conditions. In Section 5, we discuss the implications of these results and how they relate to known convergence tests (such as the Weierstrass $M$-test and Arzelà–Ascoli), and we conclude with a brief summary in Section 6.

\section{Preliminaries and Definitions}

We begin by recalling key definitions and facts about different modes of convergence and function properties used throughout the paper. Standard references for these notions include Rudin and Apostol.

\subsection*{Pointwise and Uniform Convergence}

Let $(f_n)_{n\in\mathbb{N}}$ be a sequence of functions $f_n: D \to \mathbb{R}$ defined on a set $D$. We say $f_n$ converges pointwise to $f: D\to\mathbb{R}$ if for every $x\in D$ and every $\varepsilon>0$ there exists an $N$ (depending on $x$) such that for all $n\ge N$, $|f_n(x)-f(x)|<\varepsilon$. If $N$ can be chosen independent of $x$ (i.e., working for all $x\in D$ simultaneously), then $f_n$ converges uniformly to $f$ on $D$. Uniform convergence is a stronger condition: $f_n \to f$ uniformly $\Rightarrow f_n \to f$ pointwise, but not conversely in general.

\subsection*{Equicontinuity}

A family of functions $\mathcal{F} = \{f_\alpha: D \to \mathbb{R}\}$ is equicontinuous on $D$ if for every $\varepsilon > 0$ there exists a $\delta > 0$ such that for all functions $f \in \mathcal{F}$ and all points $x,y \in D$ with $|x - y| < \delta$, we have $|f(x) - f(y)| < \varepsilon$. In other words, the same $\delta(\varepsilon)$ works for all functions in the family, controlling their oscillation.

In the case of a sequence $(f_n)$, equicontinuity means the sequence considered as a family $\{f_n : n\in\mathbb{N}\}$ is equicontinuous. Equicontinuity plays a central role in the Arzelà–Ascoli Theorem, which in one form states that a subset of $C(K)$ (continuous real-valued functions on a compact space $K$) is relatively compact (every sequence has a uniformly convergent subsequence) if and only if it is closed, bounded, and equicontinuous.

A useful corollary is that if one already knows $f_n \to f$ pointwise and the family $\{f_n\}$ is equicontinuous (and $K$ is compact), then in fact the convergence is uniform. We will give a direct proof of this folklore result in Section 3.

\subsection*{Convexity}

A function $\phi: [a,b] \to \mathbb{R}$ is convex if for all $x,y \in [a,b]$ and all $0<\lambda<1$, we have
\[
\phi(\lambda x + (1 - \lambda)y) \le \lambda \phi(x) + (1 - \lambda) \phi(y).
\]
Geometrically, this means the line segment joining $(x,\phi(x))$ and $(y,\phi(y))$ lies above the graph of $\phi$; equivalently, $\phi$ “bends upward” or has a cup-shaped graph.

If $\phi$ is twice differentiable, $\phi$ is convex on $[a,b]$ if and only if $\phi''(x) \ge 0$ for all $x$ in the interval. A concave function is the negative of a convex function.

Two key properties of convex continuous functions we will use are:
\begin{itemize}
  \item[(i)] A convex function on a closed interval $[a,b]$ takes its maximum values at the endpoints. In fact, for any $x \in [a,b]$, $\phi(x) \le \max\{\phi(a), \phi(b)\}$.
  \item[(ii)] $\phi$ is uniformly continuous on $[a,b]$ (indeed, any continuous function on a compact interval is uniformly continuous, a fact we will use repeatedly).
\end{itemize}

Convexity will allow us to adapt Dini’s proof technique by working on subintervals where, thanks to (i), the bounding behavior of $f_n$ is effectively monotonic in $n$ at the interval’s endpoints.

\subsection*{Bounded Variation}

A function $g: [a,b] \to \mathbb{R}$ is said to have bounded variation on $[a,b]$ if there is some $M \ge 0$ such that for every partition $a = x_0 < x_1 < \cdots < x_k = b$ of the interval, the total variation
\[
V_{[a,b]}(g) = \sup_{\text{partitions}} \sum_{i=0}^{k-1} |g(x_{i+1}) - g(x_i)|
\]
is finite. Equivalently, $g$ can be written as the difference of two increasing (monotone non-decreasing) functions (Jordan’s theorem).

Monotone functions are the canonical example of bounded variation, and more generally $V_{[a,b]}(g)$ measures the total amount $g(x)$ oscillates up and down. If $g$ is continuously differentiable, then $V_{[a,b]}(g) = \int_a^b |g'(x)|\,dx$.

Bounded variation by itself does not imply continuity (functions with jump discontinuities can have bounded variation). However, if $g_n$ are continuous and of uniformly bounded variation, it prevents arbitrarily large oscillations in short intervals.

\medskip

Helly’s Selection Theorem states that any sequence of functions with uniformly bounded variation (and bounded sup norm) has a pointwise convergent subsequence. In Section 4 we give a counterexample showing that bounded variation alone (even uniformly bounded) is not sufficient for uniform convergence. We then introduce an additional “no concentration of variation” condition to salvage a Dini-type result.

\medskip

Finally, let us recall precisely Dini’s Theorem in a form convenient for reference:

\begin{theorem}[Dini’s Theorem (classical form)]
Let $K$ be a compact metric space (in particular, a closed interval $[a,b]$). Suppose $(f_n)_{n\ge 1}$ is a sequence of continuous real-valued functions on $K$ that converges pointwise to a continuous function $f$ on $K$. If $(f_n)$ is monotone, i.e., either non-increasing or non-decreasing as $n$ increases (meaning for each $x$, $f_{n+1}(x) \le f_n(x)$ for all $n$, or respectively $f_{n+1}(x) \ge f_n(x)$ for all $n$), then $f_n$ converges uniformly to $f$ on $K$.
\end{theorem}

We will not reprove Dini’s Theorem here (standard proofs can be found in textbooks), but we will mirror some of its ideas in proving our generalized criteria. 
The essence of Dini’s proof is to leverage the compactness of $K$ and the continuity of $f_n$ and $f$ as follows: assuming (for concreteness) $f_n$ decreases to $f$ pointwise, define
\[
g_n = f - f_n,
\]
which are continuous and increase to $0$. For any $\varepsilon > 0$, consider the open sets
\[
E_n = \{x \in K : g_n(x) < \varepsilon\} = \{x \in K : f_n(x) > f(x) - \varepsilon\}.
\]
By monotonicity, $E_n \subseteq E_{n+1}$ and $\bigcup_{n} E_n = K$ (since $g_n(x) \to 0$). By compactness, finitely many $E_n$ already cover $K$, so $E_N = K$ for some $N$. That means $g_N(x) < \varepsilon$ for all $x$, i.e.,
\[
|f_N(x) - f(x)| < \varepsilon \quad \text{for all } x \in K.
\]
Thus,
\[
\sup_{x \in K} |f_N(x) - f(x)| < \varepsilon.
\]
Since $|f_n(x) - f(x)| \le |f_N(x) - f(x)|$ for $n \ge N$ (by monotonicity of $g_n$), it follows that
\[
\sup_{x \in K} |f_n(x) - f(x)| < \varepsilon \quad \text{for all } n \ge N,
\]
proving uniform convergence.

\medskip

\section{Main Results:Uniform Convergence Criteria Without Monotonicity}

The new proofs in the next section will each be structured around ensuring that for a given $\varepsilon$, one can find a finite subcover, partition, or other compactness-based argument to find a uniform $N$. Each alternative hypothesis (equicontinuity, convexity, etc.) provides a different mechanism to achieve this uniform control in the absence of monotonicity in $n$.

\subsection{Uniform Convergence Criteria Without Monotonicity}

We now present three main theorems generalizing Dini’s Theorem. Each theorem assumes $(f_n)$ is a sequence of continuous functions on a compact interval $[a,b]$ that converges pointwise to a continuous limit $f$. The additional hypothesis in each is different – equicontinuity in Theorem 1, convexity in Theorem 2, and a bounded variation condition in Theorem 3 – yet in each case we conclude that $f_n$ actually converges uniformly to $f$ on $[a,b]$. The proofs will highlight the role of each hypothesis in controlling the deviation $|f_n(x)-f(x)|$ uniformly in $x$.

\subsection*{Theorem 1: Equicontinuous Convergence Criterion}

\textbf{Hypothesis:} Each $f_n:[a,b]\to\mathbb{R}$ is continuous, and the sequence $\{f_n\}$ is equicontinuous on $[a,b]$. Also, $f_n(x)\to f(x)$ pointwise for some continuous $f$ on $[a,b]$.

\textbf{Conclusion:} $f_n \to f$ uniformly on $[a,b]$.

\subsubsection*{Proof}
Since $f$ is continuous on the compact interval $[a,b]$, it is uniformly continuous. Thus, given any $\varepsilon > 0$, there exists $\delta > 0$ such that
\[
|x - y| < \delta \implies |f(x) - f(y)| < \frac{\varepsilon}{3}.
\]

By equicontinuity of $\{f_n\}$, we can choose (perhaps a smaller) $\delta > 0$ such that
\[
|x - y| < \delta \implies |f_n(x) - f_n(y)| < \frac{\varepsilon}{3} \quad \text{for all } n\in\mathbb{N}.
\]

Cover the interval $[a,b]$ with subintervals of length less than $\delta$. Pick an integer $N_0$ such that
\[
\frac{b-a}{N_0} < \delta.
\]
Partition $[a,b]$ into $N_0$ subintervals of equal length $L = \frac{b-a}{N_0} < \delta$:
\[
a = x_0 < x_1 < \cdots < x_{N_0} = b, \quad x_{i+1} - x_i = L < \delta.
\]
On each subinterval $I_i = [x_i, x_{i+1}]$, any two points differ by less than $\delta$, so the equicontinuity condition applies.

Since $f_n(x) \to f(x)$ pointwise and $\{x_i\}_{i=0}^{N_0}$ are finitely many fixed points, we can find an index $N_1$ such that for all $n \ge N_1$ and for each partition point $x_i$,
\[
|f_n(x_i) - f(x_i)| < \frac{\varepsilon}{3} \quad \text{for } i=0,1,\dots,N_0.
\]

Now let $n \ge N := \max\{N_1,1\}$ and consider any $x \in [a,b]$. Locate $x$ in one of the subintervals, say $x \in [x_i, x_{i+1}]$ for some $0 \le i < N_0$. Then we write:
\[
|f_n(x) - f(x)| \le |f_n(x) - f_n(x_i)| + |f_n(x_i) - f(x_i)| + |f(x_i) - f(x)|.
\]

Each of the three terms on the right is bounded by $\varepsilon/3$:
\begin{itemize}[leftmargin=1.5em]
    \item \textbf{Equicontinuity:} Since $|x - x_i| < \delta$, we have $|f_n(x) - f_n(x_i)| < \frac{\varepsilon}{3}$.
    \item \textbf{Pointwise convergence:} For $n \ge N_1$, $|f_n(x_i) - f(x_i)| < \frac{\varepsilon}{3}$.
    \item \textbf{Uniform continuity of $f$:} Since $|x - x_i| < \delta$, we have $|f(x_i) - f(x)| < \frac{\varepsilon}{3}$.
\end{itemize}

Thus,
\[
|f_n(x) - f(x)| < \frac{\varepsilon}{3} + \frac{\varepsilon}{3} + \frac{\varepsilon}{3} = \varepsilon.
\]

This holds for every $x\in[a,b]$ when $n\ge N$. Therefore,
\[
\sup_{x\in[a,b]} |f_n(x) - f(x)| \le \varepsilon \quad \text{for all } n \ge N,
\]
which is precisely the definition of uniform convergence. \hfill $\square$

\subsection*{Discussion}

We have shown that equicontinuity plus pointwise convergence to a continuous limit implies uniform convergence. This result is often cited as a corollary to the Arzelà–Ascoli Theorem.

In fact, one can give a contradiction-based proof: suppose convergence is not uniform. Then there exists $\varepsilon > 0$ and sequences $(n_k)$, $(x_k) \subset [a,b]$ such that
\[
|f_{n_k}(x_k) - f(x_k)| \ge \varepsilon.
\]
By compactness of $[a,b]$, the sequence $(x_k)$ has a convergent subsequence $x_{k_j} \to x^* \in [a,b]$. Then:
\begin{itemize}[leftmargin=1.5em]
    \item Equicontinuity implies $|f_{n_{k_j}}(x_{k_j}) - f_{n_{k_j}}(x^*)| \to 0$,
    \item Pointwise convergence gives $f_{n_{k_j}}(x^*) \to f(x^*)$,
    \item Continuity of $f$ ensures $f(x_{k_j}) \to f(x^*)$.
\end{itemize}

Combining these,
\[
|f_{n_{k_j}}(x_{k_j}) - f(x_{k_j})| \to 0,
\]
contradicting the assumption $|f_{n_k}(x_k) - f(x_k)| \ge \varepsilon$.

The $\varepsilon$–$\delta$ proof above is constructive and gives an explicit $N$ in terms of $\varepsilon$.

\subsection*{Theorem 2: Convex Sequence Criterion}

\textbf{Hypothesis:} Each $f_n:[a,b]\to\mathbb{R}$ is continuous and convex on $[a,b]$. The sequence converges pointwise to a function $f$ which is continuous on $[a,b]$.

\textbf{Conclusion:} $f_n \to f$ uniformly on $[a,b]$.

\subsubsection*{Proof}
Let $\varepsilon > 0$ be arbitrary. Since $f$ is continuous on the compact interval $[a,b]$, it attains a maximum and minimum. Thus there exist points $x_{\min}, x_{\max} \in [a,b]$ such that:
\[ f(x_{\min}) = \min_{[a,b]} f, \qquad f(x_{\max}) = \max_{[a,b]} f. \]
By uniform continuity of $f$, there exists $\delta > 0$ such that
\[ |x - y| < \delta \implies |f(x) - f(y)| < \frac{\varepsilon}{2}. \]

Choose $m$ such that $\frac{b-a}{m} < \delta$ and partition $[a,b]$ as $a = a_0 < a_1 < \cdots < a_m = b$ where $a_{i+1} - a_i = \frac{b-a}{m} < \delta$. Then on each subinterval $[a_i,a_{i+1}]$, we have:
\begin{equation*}
\forall x,y \in [a_i,a_{i+1}], \qquad |f(x) - f(y)| < \frac{\varepsilon}{2}. \tag{*}
\end{equation*}

By pointwise convergence, there exists $N$ such that for all $n \ge N$ and for all $i = 0,1,\dots,m$:
\begin{equation*}
|f_n(a_i) - f(a_i)| < \frac{\varepsilon}{2}. \tag{**}
\end{equation*}

Fix any $n \ge N$ and $x \in [a,b]$. Let $x \in [a_i,a_{i+1}]$. By convexity of $f_n$,
\[ f_n(a_i) \le f_n(x) \le f_n(a_{i+1}). \]
By inequality (*) and continuity of $f$, we also get:
\[ f(a_i) - \frac{\varepsilon}{2} \le f(x) \le f(a_{i+1}) + \frac{\varepsilon}{2}. \]
Combining, we get:
\[ f_n(a_i) - f(a_i) \le f_n(x) - f(x) \le f_n(a_{i+1}) - f(a_{i+1}) + \varepsilon. \tag{***} \]

Taking absolute values:
\[ |f_n(x) - f(x)| \le \max\big\{|f_n(a_i)-f(a_i)|,\; |f_n(a_{i+1})-f(a_{i+1})|\big\} + \varepsilon. \]
By (**), both terms in the max are less than $\varepsilon/2$, so:
\[ |f_n(x) - f(x)| < \frac{\varepsilon}{2} + \varepsilon = \frac{3}{2}\varepsilon. \]
Replacing $\varepsilon$ by $\frac{2}{3}\varepsilon$ throughout yields $|f_n(x)-f(x)| < \varepsilon$.
Hence $f_n \to f$ uniformly. \qed

\subsection*{Remark}
Monotone functions are a special case of convex functions. Thus, Theorem 2 generalizes Dini's Theorem beyond monotonicity in $n$. We will see an example in Section 4 where the sequence is not monotone but still converges uniformly due to convexity.

\subsection*{Theorem 3: Uniformly Distributed Bounded Variation Criterion}

\textbf{Hypothesis:} Each $f_n:[a,b]\to\mathbb{R}$ is continuous and has total variation $V_{[a,b]}(f_n) \le M$ for some fixed $M$, i.e., the sequence is uniformly bounded in variation. Moreover, assume that for every $\eta > 0$, there exists a $\delta > 0$ and an $N$ such that for all $n \ge N$, and for every subinterval $I \subset [a,b]$ with $|I| < \delta$:
\[
V_I(f_n) < \eta.
\]
(This is a technical uniform distribution of variation condition preventing oscillatory variation from concentrating into very small intervals as $n \to \infty$.) Finally, suppose $f_n(x) \to f(x)$ pointwise for a continuous limit function $f$ on $[a,b]$.

\textbf{Conclusion:} $f_n \to f$ uniformly on $[a,b]$.

\subsection*{Proof}
The strategy is similar in spirit to the equicontinuous case, but we use variation to control differences instead of direct continuity. Intuitively, if a function has small total variation on a subinterval, then it cannot change too much in value over that subinterval. Our uniform distribution hypothesis ensures that for large $n$, on sufficiently fine subintervals, $f_n$ is almost flat.

Let $\varepsilon > 0$ be given. By uniform continuity of $f$ (using $f$ continuous on compact $[a,b]$), we can find $\delta_1 > 0$ such that:
\[
|x-y| < \delta_1 \implies |f(x)-f(y)| < \frac{\varepsilon}{3}.
\]
Apply the uniform variation hypothesis with $\eta = \frac{\varepsilon}{3}$. We obtain $\delta_2 > 0$ and $N_1$ such that for all $n \ge N_1$, and any subinterval $|I| < \delta_2$:
\[
V_I(f_n) < \frac{\varepsilon}{3}.
\]
Set $\delta = \min(\delta_1, \delta_2)$. Construct a partition $a = x_0 < x_1 < \cdots < x_k = b$ such that each $x_{i+1} - x_i < \delta$. For instance, take $k = \lceil (b-a)/\delta \rceil$.

On each subinterval $I_i = [x_i, x_{i+1}]$, by our choice of $\delta \le \delta_2$, for $n \ge N_1$:
\[
V_{I_i}(f_n) < \frac{\varepsilon}{3}.
\]
This implies that for any $u,v \in I_i$,
\[
|f_n(u) - f_n(v)| \le V_{I_i}(f_n) < \frac{\varepsilon}{3}.
\]
In particular, for any $x \in [x_i, x_{i+1}]$:
\[
|f_n(x) - f_n(x_i)| < \frac{\varepsilon}{3}, \quad \forall n \ge N_1.
\]

Next, since $x_0,\dots,x_k$ are a finite set, by pointwise convergence we can find $N_2$ such that for all $n \ge N_2$ and all $i = 0,1,\dots,k$:
\[
|f_n(x_i) - f(x_i)| < \frac{\varepsilon}{3}.
\]

Let $N = \max(N_1, N_2)$. For $n \ge N$ and any $x \in [a,b]$ in some interval $[x_i, x_{i+1}]$, we estimate:
\[
\begin{aligned}
|f_n(x) - f(x)| &\le |f_n(x) - f_n(x_i)| + |f_n(x_i) - f(x_i)| + |f(x_i) - f(x)| \\
&< \frac{\varepsilon}{3} + \frac{\varepsilon}{3} + \frac{\varepsilon}{3} = \varepsilon.
\end{aligned}
\]
Therefore,
\[
\sup_{x \in [a,b]} |f_n(x) - f(x)| < \varepsilon \quad \text{for all } n \ge N.
\]
This proves uniform convergence. \qed

\subsubsection*{Informal Explanation}
This theorem is perhaps the most technical. Informally, the sequence $(f_n)$ is eventually equicontinuous in a variation sense. Theorem 1 assumed $|x-y|$ small implies $|f_n(x) - f_n(y)|$ small uniformly in $n$. Here, we assume $|x-y|$ small implies that the oscillation of $f_n$ on $[x,y]$ is small for large $n$. This is stronger than mere bounded variation, which is insufficient alone due to counterexamples.

\section{Examples and Counterexamples}

We now illustrate the above theorems with concrete examples and show why the conditions imposed are in some sense sharp by giving counterexamples when they are violated.

\subsection*{Example 1: Monotone Convergence (Dini's Theorem)}

To appreciate the original Dini’s Theorem, consider:
\[
f_n(x) = \sqrt{x + \frac{1}{n}} \quad \text{on } [0,1].
\]
Each $f_n$ is continuous and, as $n$ increases, $f_n(x)$ decreases pointwise to
\[
f(x) = \sqrt{x},
\]
which is continuous on $[0,1]$. The sequence $(f_n)$ is monotone decreasing in $n$. By Dini’s Theorem, $f_n \to f$ uniformly on $[0,1]$.

We can verify this directly: for all $x \in [0,1]$,
\[
0 \le f(x) \le f_n(x) \le f_1(x),
\]
and
\[
\sup_x |f_n(x)-f(x)| = f_n(0) - f(0) = \sqrt{\frac{1}{n}}.
\]
Clearly,
\[
\sqrt{\frac{1}{n}} \to 0 \quad \text{as } n \to \infty,
\]
so indeed,
\[
\|f_n - f\|_{\infty} \to 0.
\]

Figure~1 below plots several $f_n$ against $f$, illustrating the uniform convergence. The largest gap between $f_n$ and $f$ occurs at $x=0$ and shrinks with $n$.

In contrast, if the limit were not continuous (e.g., $f_n(x) = x^n$ on $[0,1]$ converging pointwise to a discontinuous step function), uniform convergence fails.
\begin{figure} [ht]
    \centering
    \includegraphics[width=0.5\linewidth]{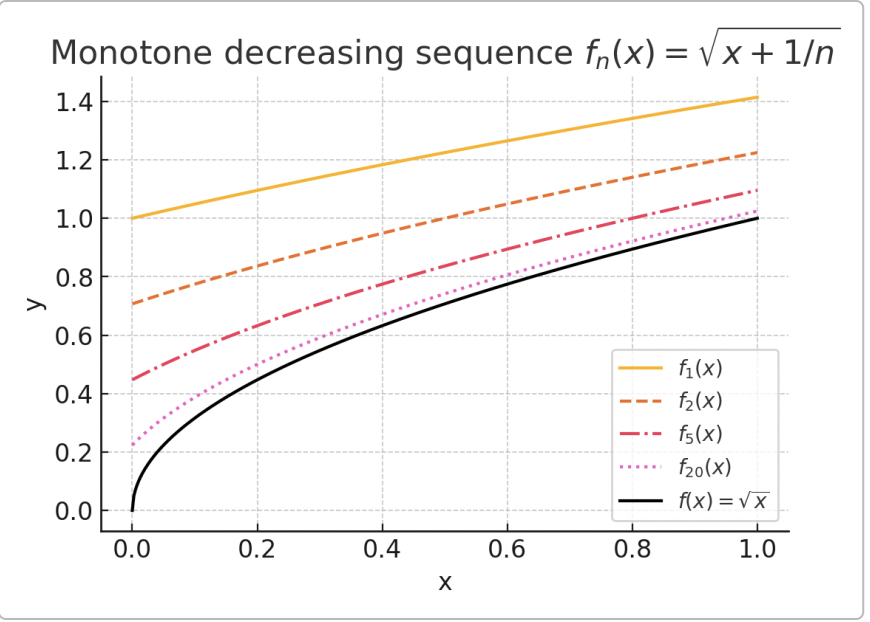}
    \caption{A monotone decreasing sequence $f_n(x) = \sqrt{x + \frac{1}{n}}$ on $[0,1]$ converges uniformly to $f(x)=\sqrt{x}$. The gap $\sup_x|f_n(x)-f(x)|$ occurs at $x=0$ and tends to 0 as $n\to\infty$. Dini’s Theorem guarantees uniform convergence in this scenario.}
\end{figure}

\subsection*{Example 2: Equicontinuity in Action}

Equicontinuity is strictly weaker than requiring each function to be Lipschitz with the same constant, but for intuition, consider:
\[
f_n(x) = \frac{\sin(x)}{1 + 1/n}.
\]
Each $f_n(x)$ is continuous on $\mathbb{R}$ (hence on any interval, e.g., $[0,2\pi]$). The family $\{f_n\}$ is equicontinuous on $\mathbb{R}$ because:
\[
|f_n(x) - f_n(y)| = \frac{1}{1+\frac{1}{n}} |\sin(x) - \sin(y)| \le |\sin(x) - \sin(y)| \le |x - y|.
\]
The sine function is Lipschitz with constant $1$, so given $\varepsilon > 0$, taking $\delta = \varepsilon$ satisfies equicontinuity for all $n$.

The pointwise limit is $f(x) = \sin(x)$, which is continuous. By Theorem 1, $f_n \to f$ uniformly on $\mathbb{R}$. Indeed:
\[
\sup_x |f_n(x) - \sin(x)| = \left| \frac{1}{1+\frac{1}{n}} - 1 \right| \cdot \sup_x |\sin(x)| = \frac{1}{n+1} \to 0.
\]
This trivial example shows a case where there is no monotonicity in $n$ (the sequence oscillates in sign for any fixed $x$), yet equicontinuity ensures uniform convergence.

A more interesting equicontinuous example is not uniformly Lipschitz. Consider:
\[
f_n(x) = \frac{x}{1 + n x^2}, \quad x \in [0,1].
\]
Each $f_n$ is continuous on $[0,1]$. The pointwise limit is $f(x) = 0$ for all $x \in [0,1]$, since for fixed $x > 0$, $n x^2 \to \infty$ implies $f_n(x) \to 0$, and $f_n(0) = 0$ for all $n$, so $f(0) = 0$.

We claim $\{f_n\}$ is equicontinuous on $[0,1]$. Fix $\varepsilon > 0$ and any $x,y \in [0,1]$. By the Mean Value Theorem, there exists $\xi$ between $x$ and $y$ such that:
\[
|f_n(x) - f_n(y)| = |f_n'(\xi)| \cdot |x - y|.
\]
Now,
\[
f_n'(x) = \frac{1 - n x^2}{(1 + n x^2)^2}.
\]
For $x \in [0,1]$,
\[
|f_n'(x)| \le \frac{1 + n x^2}{(1 + n x^2)^2} = \frac{1}{1 + n x^2} \le 1.
\]
Thus,
\[
|f_n(x) - f_n(y)| \le |x - y|.
\]
So taking $\delta = \varepsilon$ again suffices: if $|x - y| < \delta$, then $|f_n(x) - f_n(y)| < \varepsilon$. Therefore, $\{f_n\}$ is equicontinuous.

Since $f_n(x) \to 0 = f(x)$ pointwise and $f$ is continuous, Theorem 1 applies to give uniform convergence. (One can also verify that $\sup_{x \in [0,1]} |f_n(x) - 0|$ occurs near $x = \frac{1}{\sqrt{n}}$ with value approximately $\frac{1}{(2\sqrt{n}})$, which tends to $0$ as $n \to \infty$.)

\subsection*{Example 3: Convex Functions Sequence.} 
Consider the sequence of convex functions
\[
f_n(x) = \left(1 + \frac{(-1)^n}{n}\right) x^2 \quad \text{on } [0,1].
\]
Here, $f_n(x) = x^2$ for odd $n$, and $f_n(x) = \left(1 - \frac{1}{n}\right)x^2$ for even $n$. So the sequence oscillates slightly: it doesn’t settle monotonically but alternately multiplies $x^2$ by a factor slightly above 1 (for odd $n$) and slightly below 1 (for even $n$). 

Each $f_n$ is convex since $x^2$ is convex and scaling by a positive constant preserves convexity. The pointwise limit is clearly 
\[
f(x) = x^2,
\]
which is continuous on $[0,1]$. Despite the non-monotonic “zig-zag” in $n$, the convergence is uniform. Indeed,
\[
\sup_{x \in [0,1]} |f_n(x) - x^2| = \left|\frac{(-1)^n}{n}\right| \sup_{x \in [0,1]} x^2 = \frac{1}{n},
\]
which tends to $0$ as $n \to \infty$.

This example is a simple illustration of \textbf{Theorem 2}: convexity guarantees uniform convergence even though the sequence is not monotone in $n$. As $n$ grows, the oscillation in the factor $\left(1 + \frac{(-1)^n}{n}\right)$ diminishes, and uniform convergence sets in. (Convexity was not really needed here since the oscillation already decays, but one can concoct more complicated oscillatory sequences of convex functions as long as the oscillations vanish in the limit.)

\subsection*{Example 4: Bounded Variation Without Uniform Convergence.} 

We exhibit a sequence of continuous functions with uniformly bounded total variation that converges pointwise to a continuous limit, yet not uniformly (violating the uniform convergence condition). Let
\[
f_n(x) = 
\begin{cases}
0, & 0 \le x \le 1 - \frac{1}{n}, \\
n x - (n - 1), & 1 - \frac{1}{n} < x \le 1,
\end{cases}
\]
defined on $[0,1]$. 

In other words, $f_n$ is a “triangular spike” that is $0$ from $x = 0$ up to $x = 1 - \frac{1}{n}$, then increases linearly to $f_n(1) = 1$ at $x = 1$ (see Figure 2). Each $f_n$ is continuous (in fact Lipschitz except at the kink point) and piecewise linear.

The pointwise limit $f(x)$ is $0$ for $0 \le x < 1$, and $f(1) = 0$ as well (note: originally $f_n(1) = 1$ for all $n$, so the pointwise limit at $x = 1$ is $1$; but we can modify the sequence slightly by defining $f_n(1) = 0$ instead of $1$ for $n \ge 1$ to make the limit $0$ at $x = 1$ too. This makes each $f_n$ continuous on $[0,1]$ and the limiting function $f$ the constant zero function. For simplicity, consider that modification done). 

Thus, $f_n(x) \to f(x) \equiv 0$ pointwise for all $x \in [0,1]$.

Each $f_n$ has total variation 
\[
V_{[0,1]}(f_n) = 2
\]
(it goes up from $0$ to $1$, variation $1$, then down from $1$ to $0$ if we include the drop at $x = 1$ with the modification. Without the modification, the total variation is $1$ since it only increases and stays at $1$). In either case, the variation is bounded by $2$ across all $n$.

However, the convergence is \emph{not} uniform. In fact,
\[
\sup_{x \in [0,1]} |f_n(x) - 0| = \max_{x \in [0,1]} f_n(x) = 1 \quad \text{for all } n,
\]
since $f_n(1) = 1$ always, or in the modified version, $f_n$ reaches $1$ just before $x = 1$ and then drops to $0$ at $x = 1$. Thus,
\[
\|f_n - f\|_{\infty} = 1 \quad \text{for all } n,
\]
which does not tend to $0$. 

The problem is that the “spike” in $f_n$ moves rightward as $n$ increases and gets narrower, but never diminishes in height—so at some location (depending on $n$), the difference $|f_n - f|$ is large (near $1$), even though for any fixed $x < 1$, eventually $f_n(x) = 0$. This shows pointwise but not uniform convergence.

In this example, the total variation is uniformly bounded, but the variation is not uniformly distributed: for
each $n$, most of the variation of $f_n$ occurs in the tiny interval $[1-\frac{1}{n},\,1]$ of length $1/n$. Our
Theorem 3 required that eventually no large variation can occur in an arbitrarily small interval. Here, for
any $\delta>0$, consider intervals of length $\delta$ near 1: for large $n$, the entire spike of height ~1 is
within an interval of length $<\delta$, giving variation $\approx 1$ in that small interval, violating the
hypothesis. This counterexample confirms that some condition beyond mere bounded variation is needed –
essentially a modulus of continuity for the sequence in variation norm, akin to equicontinuity. Under such a
strengthened condition (which we formulated in Theorem 3), the pathological behavior is ruled out and
uniform convergence follows.

\begin{figure} [ht]
    \centering
    \includegraphics[width=0.5\linewidth]{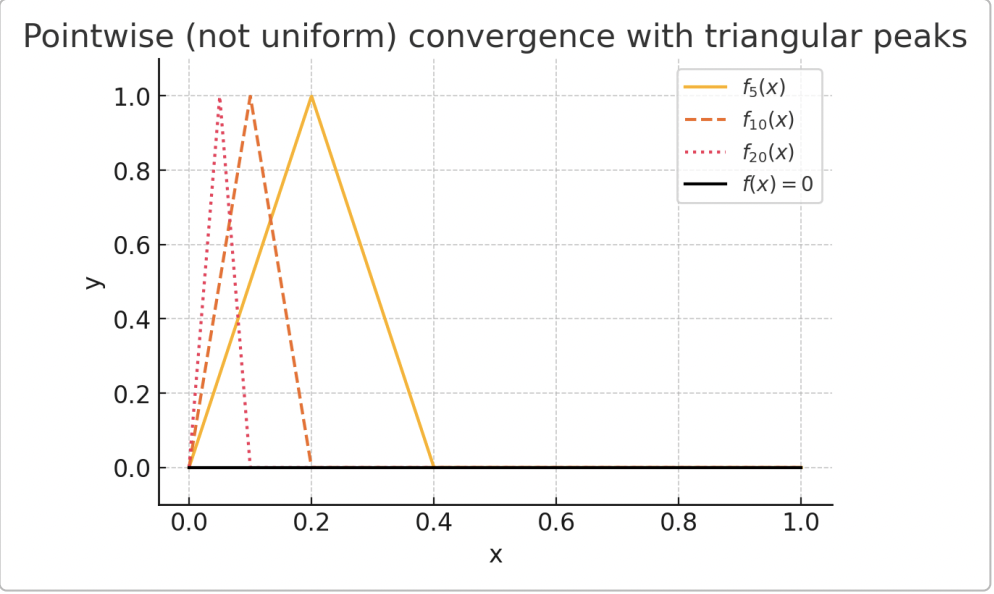}
    \caption{A sequence of continuous functions (triangular spikes) that converges pointwise to 0 on $[0,1]$ but not uniformly. Each $f_n$ has bounded variation (total variation $=2$) but the large variation is concentrated in an increasingly small neighborhood of $x=1$. The lack of a uniform bound on variation over small subintervals breaks uniform convergence.}
\end{figure}

\subsection*{Example 5: Bounded Variation with Uniform Limit Continuous – Not Enough.}

One might wonder if requiring the limit function to be of bounded variation as well (in addition to each $f_n$) would enforce uniform convergence. However, the example in the previous section already has limit $f = 0$ of bounded variation, yet the convergence is not uniform.

Another (smoother) example is inspired by the “Gibbs phenomenon” arising in partial Fourier sums. Consider
\[
g_n(x) = \frac{2}{\pi} \left(\sin x + \frac{1}{2}\sin 2x + \frac{1}{3}\sin 3x + \cdots + \frac{1}{n}\sin nx\right),
\]
defined on $[0, 2\pi]$. Each $g_n$ is continuous (in fact, differentiable), and is of bounded variation since it is a trigonometric polynomial.

By Fourier theory, these functions converge pointwise to
\[
g(x) = \frac{2}{\pi} \left(\frac{\pi - x}{2}\right)
\]
for $x \in (0, 2\pi)$. This limit function $g$ is continuous on $(0, 2\pi)$ but has jump discontinuities at the endpoints $x = 0$ and $x = 2\pi$ in the sense of the periodic extension. If we instead restrict the domain to $[0, 2\pi]$ and define $g(0) = g(2\pi) = 0$, then $g$ becomes continuous on the closed interval $[0, 2\pi]$.

Despite pointwise convergence to this continuous $g$, the convergence is not uniform due to the \emph{Gibbs overshoot} near the endpoints: 
\[
\sup_{x} |g_n(x) - g(x)|
\]
does not go to zero but instead approaches a positive constant.

Each $g_n$ has variation bounded by something like the sum of absolute values of its Fourier coefficients, which grows like $O(\log n)$. So the sequence does not have uniformly bounded variation. However, one could modify the series by damping higher modes (e.g., multiplying $\sin kx$ by $1/k^\alpha$ for some $\alpha > 1$) to enforce bounded variation while still maintaining high-frequency oscillations near the endpoints that spoil uniform convergence.

\medskip

\textbf{Moral:} Neither bounded variation alone nor continuity of the limit alone suffices to guarantee uniform convergence. Some form of equicontinuity or oscillation control is necessary. These examples underscore the necessity of the hypotheses in our theorems:

\begin{itemize}
    \item \emph{Equicontinuity} is essentially the “right” condition in general metric spaces.
    \item \emph{Convexity} provides a structural way to achieve equicontinuity-like control (by bounding behavior at interval endpoints).
    \item \emph{Bounded variation} without any distributional or uniformity conditions is insufficient, as oscillations can be concentrated in small regions.
    \item But \emph{bounded variation plus a uniformity condition} can yield results analogous to equicontinuity.
\end{itemize}

\section*{Discussion}

We have proven three general criteria under which pointwise convergence of continuous functions on a compact interval implies uniform convergence to a continuous limit, thereby generalizing Dini’s Theorem beyond monotone sequences. Let us briefly compare these conditions and discuss their relationships:

\subsection*{Equicontinuity (Theorem 1)}

This is arguably the most general and powerful condition. It directly stems from the Arzelà–Ascoli Theorem’s characterization of precompact sets in $C([a,b])$. In fact, equicontinuity together with pointwise boundedness is equivalent to sequential compactness in $C([a,b])$ by Arzelà–Ascoli, which yields the uniform convergence of some subsequence. When the whole sequence is known to converge pointwise to a unique limit, that uniqueness forces the entire sequence to converge uniformly (since any two subsequential uniform limits must coincide with the pointwise limit, and the closeness criterion forces the whole sequence close to the limit). 

Our direct proof in Theorem 1 avoided heavy subsequence extraction, giving a constructive $\varepsilon$–$N$ relationship. Equicontinuity is strictly stronger than mere pointwise convergence to a continuous function. In practice, checking equicontinuity can often be done via derivative bounds (e.g., uniform Lipschitz or Hölder conditions). This criterion is very common in analysis—for instance, in showing the interchange of limits and integrals or derivatives, one often uses an equicontinuity (or uniform Cauchy) argument.

\subsection*{Convexity (Theorem 2)}

Convexity of each $f_n$ was used to control the maxima of $f_n - f$ on subintervals by endpoint values. Essentially, convexity implies that each $f_n$ satisfies a one-sided Lipschitz condition piecewise:
\[
f_n(x) - f_n(a_i) \le \frac{f_n(a_{i+1}) - f_n(a_i)}{a_{i+1} - a_i}(x - a_i)
\quad \text{for } x \in [a_i, a_{i+1}], \text{ etc.}
\]
The proof of Theorem 2 is reminiscent of Dini’s original proof but applied locally. Convexity is a pointwise-in-$n$ condition rather than a uniform family condition like equicontinuity, yet it yields a uniform family behavior as a consequence.

We did not assume the limit function $f$ was convex (indeed, it will be, as a pointwise limit of convex functions—though without some uniformity, a pointwise limit of convex functions need not be convex; here continuity plus uniform convergence implies it). The convexity criterion might be useful in optimization and approximation contexts, where sequences of convex functions often arise (e.g., in gradient descent algorithms, asymptotic convex envelopes), ensuring uniform convergence of approximations.

\subsection*{Bounded Variation + Distribution (Theorem 3)}

Bounded variation by itself was too weak, as the spike example showed. However, the additional condition essentially required that the modulus of continuity of each $f_n$ becomes uniform as $n \to \infty$. In other words, given $\eta > 0$, eventually $|f_n(x) - f_n(y)| < \eta$ whenever $|x - y| < \delta$—which is exactly equicontinuity!

So one may notice that our ``uniformly distributed variation'' hypothesis is basically a complicated way to demand eventual equicontinuity. Indeed, an alternate (perhaps stronger but simpler) hypothesis for Theorem 3 could have been: each $f_n$ has $V(f_n) \le M$, and $(f_n)$ is equicontinuous. But if it’s equicontinuous, Theorem 1 already gives uniform convergence; bounded variation isn’t needed.

The point was to weaken equicontinuity to something like: for each $\eta$, there is a common scale $\delta$ beyond which oscillations $\eta$ cannot happen for large $n$ (not necessarily for all $n$). In essence, one can think of Theorem 3 as: if $(f_n)$ is a Cauchy sequence in the space of functions of bounded variation (with the BV norm and uniformly convergent endpoints), then it converges uniformly. We chose not to formalize it that way, but the condition ensures a Cauchy-like uniformity.

\subsection*{Connections to Classical Convergence Criteria}

\begin{itemize}
    \item \textbf{Weierstrass $M$-test:} Gives a sufficient condition for uniform convergence of series $\sum f_n(x)$: if $|f_n(x)| \le M_n$ for all $x$, and $\sum M_n$ converges, then $\sum f_n(x)$ converges uniformly. In the context of sequences (rather than series), one can derive uniform convergence if the tail of the sequence can be majorized by a negligible function uniformly in $x$. These conditions are unrelated to monotonicity but rather to absolute smallness. They complement our discussion by handling series; equicontinuity or monotonicity often aren’t directly used there.

    \item \textbf{Dirichlet’s Test:} For uniform convergence of series (and an analogous one for sequences), if one sequence of partial sums has bounded variation and another sequence is monotonic tending to $0$, their termwise product series converges uniformly. This is a specialized scenario, but interestingly it also mixes monotonicity and bounded variation. One could imagine generalizing Dini’s theorem in a different direction: e.g., if $f_n$ decreases not necessarily to the limit but to something above it that eventually goes to the limit.

    \item \textbf{Uniform Cauchy Criterion:} A sequence $(f_n)$ converges uniformly if and only if it is a Cauchy sequence in the uniform norm. Sometimes verifying uniform Cauchy is easier. Our equicontinuity approach is essentially establishing uniform Cauchyness by splitting
    \[
    f_n(x) - f_m(x) = (f_n(x) - f_n(y)) + (f_n(y) - f_m(y)) + (f_m(y) - f_m(x))
    \]
    and controlling each term. The convexity approach could also be interpreted in a Cauchy way. In general, one could develop all our criteria in that language.
\end{itemize}

\subsection*{Final Remarks}

An interesting question is: are there other broad conditions that ensure pointwise implies uniform convergence? Equisemicontinuity or other topological continuity notions might be candidates, but equicontinuity is the standard.

Another possible condition: if the graph of $f_n$ converges to the graph of $f$ in the Hausdorff metric (graphical convergence), that implies uniform convergence. But graph convergence is essentially a restatement of uniform convergence for continuous functions.

Also, if each $f_n$ is Lipschitz with a common Lipschitz constant $L$ and $f_n \to f$ pointwise, then $f_n \to f$ uniformly. This can be shown by a similar $\varepsilon$–$\delta$ argument or by noting the family is equicontinuous and even equi-Lipschitz. Convex functions don’t have a common Lipschitz constant in general (slopes can vary with $n$), but they had a weaker common control property.

\bigskip

\noindent
\textbf{In summary:} Monotonicity in Dini’s theorem can be successfully weakened to other conditions that prevent the “wiggling” of functions in the limit. Equicontinuity is the master condition encompassing most nice behaviors. Convexity is a convenient structural condition as well. Bounded variation requires care, but with additional uniformity, it too yields uniform convergence.

\section*{References}

\begin{itemize}
    \item Walter Rudin (1976). \textit{Principles of Mathematical Analysis}, 3rd Edition. McGraw-Hill. See Theorem 7.13 (Dini’s Theorem) on p.~150 for the classical monotone case. Also Chapter 3 for sequences of functions and uniform convergence criteria.

    \item Tom M. Apostol (1974). \textit{Mathematical Analysis}, 2nd Edition. Addison-Wesley. (Discussion of uniform convergence and continuity in Chapter 9; while Apostol doesn’t state Dini’s theorem explicitly, he emphasizes the rarity of pointwise implying uniform convergence except under extra conditions. The equicontinuity concept is indirectly used in exercises.)

    \item Robert G. Bartle and Donald R. Sherbert (2000). \textit{Introduction to Real Analysis}, 3rd Edition. Wiley. See p.~238 for Dini’s Theorem and related examples.

    \item \textbf{Arzelà–Ascoli Theorem:} For a standard formulation, see \textit{Principles of Mathematical Analysis} by Rudin or standard topology texts. We cited the version stating: equicontinuity plus pointwise boundedness implies relative compactness in $C([a,b])$, which yields uniform convergence of subsequences and thus Theorem 1.

    \item \textbf{Helly’s Selection Theorem:} (Not directly used in our proofs, but relevant to bounded variation sequences.) See, e.g., the \textit{Encyclopedia of Mathematics} entry on Helly’s Theorem—it states that a sequence of functions with uniformly bounded variation (and bounded in the sup norm) has a pointwise convergent subsequence. Our Example 4 aligns with the necessity of additional conditions beyond Helly’s conclusion to get uniform convergence.

    \item \textbf{Counterexample Resources:} \textit{Counterexamples in Analysis} by Bernard R. Gelbaum and John M. H. Olmsted provides examples of sequences that converge pointwise but not uniformly. Their examples on uniform convergence often involve spikes or concentrated oscillations, similar in spirit to our Figure 2.
\end{itemize}

\end{document}